\magnification=1200
\input amstex
\UseAMSsymbols

\input miniltx
 \expandafter\def\expandafter\+\expandafter{\+}
 \input url.sty
 
\centerline{\bf  ROOTS AND LOGARITHMS OF MULTIPLIERS}

\bigskip

\centerline{\bf  Jingbo Xia, Congquan Yan, Danjun Zhao and Jingming Zhu}

{\input amsppt.sty
\topmatter
\thanks
\endthanks
\endtopmatter
}

\noindent
{\bf Abstract.}   By now it is a well-known fact that if  $f$  is a multiplier for the Drury-Arveson space  
$H^2_n$,  and if there is a  $c > 0$  such that   $|f(z)| \geq c$  for every  $z \in {\bold B}$,  then the reciprocal 
function  $1/f$  is also a multiplier for   $H^2_n$.  We show that for such an  $f$  and for every  $t \in {\bold R}$,  
$f^t$  is also a multiplier for   $H^2_n$.  We do so by deriving a differentiation formula for   $R^m(f^th)$. 
Moreover, by this formula the same result holds for spaces  ${\Cal H}_{m,s}$ of the Besov-Dirichlet type.    
The same technique also gives us the result that for a non-vanishing multiplier  $f$  of  $H^2_n$,   
$\log f$  is a multiplier of   $H^2_n$  if and only if  $\log f$  is bounded on  ${\bold B}$.
\bigskip
\centerline{\bf  1.  Introduction}
\medskip
Recall that the Drury-Arveson space
$H^2_n$   is the Hilbert space of analytic functions on the unit ball    
${\bold B} = \{z \in {\bold C}^n : |z| < 1\}$  that has the function  
$$
K_z(\zeta )  =  {1\over 1 - \langle \zeta ,z\rangle }
$$
\noindent
as its reproducing kernel [1,4,9].   A function  $f \in H^2_n$  is said to be a {\it multiplier} for  $H^2_n$  
if it has the property that  $fh \in H^2_n$  for every  $h \in H^2_n$.  Let  ${\Cal M}$  denote the collection of the 
multipliers for the Drury-Arveson space  $H^2_n$.  It is well known that for each  $f \in {\Cal M}$,  the multiplication operator  
$M_fh = fh$,  $h \in H^2_n$,  is necessarily bounded on  $H^2_n$  [1].
\medskip
In a remarkable paper [3], Costea, Sawyer and Wick proved the corona theorem for  ${\Cal M}$.  This 
in particular implies the one-function corona theorem:  
\medskip
\noindent
{\bf Theorem 1.1.}  {\it If} $f$  {\it is a multiplier for}   
$H^2_n$,  {\it and if there is a}  $c > 0$  {\it such that}   $|f(z)| \geq c$  {\it for every}  $z \in {\bold B}$,  {\it then the reciprocal function}  
$1/f$  {\it is also a multiplier for}   $H^2_n$. 
\medskip
Shortly after [3] was published, interest arose to find a ``short" proof of Theorem 1.1,  a proof that 
does not rely on the heavy analytical tools in [3].  The first proof that meets the criterion of being ``short" was given in [5].   
While elementary, the proof in [5] of the  one-function corona theorem still required 
one analytical tool, namely the von Neumann inequality for commuting tuples of row contractions.   Later a second ``short" proof of 
Theorem 1.1 was given in [10], but it also required an analytical tool in the form of  [10, Lemma 5.3],  
which is similar in spirit to the von Neumann inequality for commuting row contractions.  
Also see [6, Corollary 5.1]  and the discussion on this subject in [7].
\medskip
 Then in [2],  Cao, He and Zhu gave a surprising proof of Theorem 1.1 by a formula.
In addition to the elegance of proving a result by a formula,
their approach has the advantage that it also works on spaces where the above-described analytical tools are not available.  
The main ingredient in [2] is the radial derivative.  
\medskip
Throughout the paper, we denote
$$
R  =  z_1\partial _1 + \cdots + z_n\partial _n.
$$
\noindent
In other words,  $R$  denotes the usual radial derivative on  ${\bold B}$.
\medskip
Let   $m \in {\bold N}$  and   $-1 < s < \infty $.  We define  ${\Cal H}_{m,s}$  to be the the collection of analytic functions  
$h$  on  ${\bold B}$  satisfying the condition  $\|h\|_{m,s} < \infty $,  where the norm  $\|\cdot \|_{m,s}$  is defined by the 
formula
$$
\|h\|_{m,s}^2  =  |h(0)|^2  +  \int _{\bold B}|(R^mh)(z)|^2(1 - |z|^2)^sdv(z).
$$
\noindent
Then   ${\Cal H}_{m,s}$  is a Hilbert space of analytic functions on  ${\bold B}$.   Obviously,  
${\Cal H}_{m,s}$   is a multi-variable, higher-derivative analogue of the Dirichlet space.  See [12].
\medskip
A function  $f \in {\Cal H}_{m,s}$  is said to be a {\it multiplier}  for     ${\Cal H}_{m,s}$  if   $fh \in   {\Cal H}_{m,s}$  
for every   $h \in   {\Cal H}_{m,s}$.  If  $f$  is a multiplier for  ${\Cal H}_{m,s}$,
then it follows easily from the closed-graph theorem that the multiplication operator  $M_fh = fh$  is necessarily bounded 
on  ${\Cal H}_{m,s}$.  The operator norm  $\|M_f\|$  on  ${\Cal H}_{m,s}$   is generally called the multiplier norm of  $f$.  
\medskip
We recall the 
following result of Cao, He and Zhu:
\medskip
\noindent
{\bf Theorem 1.2.}  [2]  {\it Let }  $m \in {\bold N}$  {\it and}  $-1 < s < \infty $,  {\it and let}  $f$  {\it be a multiplier for}    ${\Cal H}_{m,s}$.  
{\it If there is a}  $c > 0$  {\it such that}  $|f(z)| \geq c$  {\it for every}  $z \in {\bold B}$,  {\it then}    $1/f$  {\it is also a multiplier for}   ${\Cal H}_{m,s}$. 
\medskip
Given any  $n \in {\bold N}$,  let  $m_0 \in {\bold N}$  and  $k_0 \in {\bold Z}_+$  be such that  $2m_0 - k_0 = n$.  
Then by a well-known estimate,  there are constants  $0 < c \leq C < \infty $  such that  
$$
c\|h\|_{m_0,k_0}  \leq  \|h\|_{H^2_n}  \leq  C\|h\|_{m_0,k_0}
\tag 1.1                  
$$
\noindent
for every   $h \in H^2_n$.  Therefore Theorem 1.2 implies Theorem 1.1. 
In [2], Theorem 1.2 was proved by using the following formula:
$$
R^m(f^{-1}h)  =  \sum _{k=0}^m(-1)^{m-k}{(m+1)!\over k!(m+1-k)!}{1\over f^{m-k+1}}R^m(f^{m-k}h)
\tag 1.2
$$
\noindent
if  $f$,  $h$  are analytic on  ${\bold B}$  and if  $f$  does not vanish on  ${\bold B}$.
\medskip
Given Theorems 1.1 and 1.2,  a natural thought that comes to mind is, what about the roots, or fractional powers, of  $f$?  
For the unit ball,  it is well know that if  $f$  is a non-vanishing analytic function on  ${\bold B}$,  then there is an 
analytic function  $g$  on  ${\bold B}$  such that  $f = e^g$,  i.e.,  $f$  has a logarithm.  
Thus if $f$  is a non-vanishing analytic function on  ${\bold B}$,  
then for every real number  $t \in {\bold R}$,  we have a well-defined analytic function  $f^t$  on  ${\bold B}$.
\medskip
Now, if  $f$  is a multiplier for  $H^2_n$,  and if there is a  $c > 0$  such that  $|f(z)| \geq c$  for every  $z \in {\bold B}$,  
then for a general  $t \in {\bold R}$,  is  $f^t$  a multiplier for  $H^2_n$?  Furthermore,  
what about the multipliers of  ${\Cal H}_{m,s}$?    We have the following results:
\medskip
\noindent
{\bf Theorem 1.3.}   {\it Let }  $m \in {\bold N}$  {\it and}  $-1 < s < \infty $.  {\it Let}  $f$  {\it be a multiplier for}    ${\Cal H}_{m,s}$,    
{\it and suppose that there is a}  $c > 0$  {\it such that}  $|f(z)| \geq c$  {\it for every}  $z \in {\bold B}$. 
{\it Then for every}   $t \in {\bold R}$,    {\it the function}  $f^t$  {\it is also a multiplier for}   ${\Cal H}_{m,s}$. 
\medskip
It is also natural to consider the logarithm of a non-vanishing multiplier:
\medskip
\noindent
{\bf Theorem 1.4.}   {\it Let }  $m \in {\bold N}$  {\it and}  $-1 < s < \infty $.  {\it Let}  $f$  {\it be a multiplier for}    ${\Cal H}_{m,s}$,    
{\it and suppose that} $f$  {\it does not vanish on}  ${\bold B}$. 
{\it Then}   $\log f$    {\it is a multiplier for}   ${\Cal H}_{m,s}$   {\it if and only if}   $\log f$  {\it is bounded on}  ${\bold B}$. 
\medskip
By (1.1), Theorem 1.3 immediately implies that if $f$  is a multiplier for  the Drury-Arveson space
$H^2_n$  and if there is a $c > 0$  such that
$|f(z)| \geq c$  for every  $z \in {\bold B}$,  then for every  $t \in {\bold R}$,    the function $f^t$   is also a multiplier for
$H^2_n$.  Similarly, it follows from (1.1) and Theorem 1.4  that if   $f$  is a non-vanishing multiplier for the 
Drury-Arveson space  $H^2_n$,  then  $\log f$  is a multiplier for  $H^2_n$  if and only if   
$\log f$  is bounded on  ${\bold B}$.
\medskip
Note that if  $f$  is a multiplier,  then the sequence of functions  $f, f^2, \dots ,f^k, \dots $  are all multipliers.  
Our proofs of Theorems 1.3 and 1.4  are based on two differentiation formulas,  
which are the main results of the paper:
\medskip
\noindent
{\bf Proposition 1.5.}    {\it Given any}  $m \in {\bold N}$  {\it and}    $t \in {\bold R}$,  {\it there are real numbers}   
$\rho _1$, $\dots $, $\rho _m$  
{\it such that for all analytic functions}  $f$, $h$  {\it on}  ${\bold B}$,  {\it  if}  $f$  {\it does not vanish on}  ${\bold B}$,  {\it then}
$$
R^m(f^th) - f^tR^mh = f^t\sum _{k=1}^m{\rho _k\over f^k}\{R^m(f^kh) - f^kR^mh\}.
\tag 1.3
$$
\medskip
\noindent
{\bf Proposition 1.6.}    {\it Given any}  $m \in {\bold N}$,  {\it there are rational numbers}   $a_1$, $\dots $, $a_m$  
{\it such that for all analytic functions}  $f$, $h$  {\it on}  ${\bold B}$,  {\it  if}  $f$  {\it does not vanish on}  ${\bold B}$,  {\it then}
$$
R^m((\log f)h) - (\log f)R^mh = \sum _{k=1}^m{a_k\over f^k}\{R^m(f^kh) - f^kR^mh\}.
\tag 1.4
$$
\medskip
In Sections 3 and 4 we will see that the coefficients  $\rho _1$, $\dots $, $\rho _m$    and   $a_1$, $\dots $, $a_m$  
in Propositions 1.5 and 1.6 are explicitly determined.

\bigskip
\centerline{\bf  2.  Proofs of Theorems 1.3 and 1.4}
\medskip
It is well known that if  $f$  is a multiplier for the  Drury-Arveson space  $H^2_n$,  then  $f$  is necessarily bounded 
on the unit ball  ${\bold B}$  [1],  i.e.,  $\|f\|_\infty < \infty $.   
The same is true for the multipliers of the spaces  ${\Cal H}_{m,s}$  introduced in Section 1,  
although this fact is not as well known.  For completeness, we first establish
\medskip
\noindent
{\bf Lemma 2.1.}   {\it Let}  $m \in {\bold N}$  {\it and}  $-1 < s < \infty $.  {\it If}  $f$  {\it is a multiplier for}  ${\Cal H}_{m,s}$,  
{\it then}   $\|f\|_\infty < \infty $.
\medskip
\noindent
{\it Proof}.   Write  $\|\cdot \|_s$  for the norm of the weighted Bergman space  $L^2_a({\bold B},(1-|z|^2)^sdv(z))$.   
By the definition of the norm  $\|\cdot \|_{m,s}$,  there is a  $0 < C_s < \infty $  such that   
$\|h\|_s \leq C_s\|h\|_{m,s}$  for every   $h \in {\Cal H}_{m,s}$.  Therefore for each given  $w \in {\bold B}$,  the map  
$$
h  \mapsto h(w)  
$$
is a bounded linear functional on   ${\Cal H}_{m,s}$.  Hence for each  $w \in {\bold B}$,  there is a  $K^{(m,s)}_w \in {\Cal H}_{m,s}$  
such that  
$$
h(w)  =  \langle h,K^{(m,s)}_w\rangle _{m,s}   \quad   \text{for every}  \ \   h \in {\Cal H}_{m,s}.
$$
\noindent
Let  $f$  be a multiplier for   ${\Cal H}_{m,s}$.  For every pair of  $w, z \in {\bold B}$,  we have
$$
\align
\langle M_f^\ast K^{(m,s)}_w,K^{(m,s)}_z\rangle _{m,s}  &=  \langle K^{(m,s)}_w,M_fK^{(m,s)}_z\rangle _{m,s}  
=  \overline{f(w)K^{(m,s)}_z(w)}  \\
&=   \overline{f(w)}\langle K^{(m,s)}_w,K^{(m,s)}_z\rangle _{m,s}.
\endalign
$$
\noindent
From this we deduce   $M^\ast _fK^{(m,s)}_w = \overline{f(w)}K^{(m,s)}_w$  for every   $w \in {\bold B}$.   Thus  
$\|f\|_\infty \leq \|M^\ast _f\| = \|M_f\| < \infty  $. This completes the proof.   $\square $
\medskip
\noindent
{\it Proof of Theorem} 1.3.  
Let  $f$  be a multiplier for  ${\Cal H}_{m,s}$, and suppose that there is a  $c > 0$  such that
$|f(z)| \geq c$  for every  $z \in {\bold B}$.  Then   $\|f^{t-k}\|_\infty  < \infty $, whether  $t - k \geq 0 $  or  $t - k < 0$.
\medskip  
Applying Proposition 1.5   to  $f$  and to each  $h \in {\Cal H}_{m,s}$,  we obtain  
$$
R^m(f^th)  =  \sum _{k=0}^m\rho _kf^{t-k}R^m(f^kh),
$$
\noindent
where  $\rho _0 = 1 - \rho _1 - \cdots - \rho _m$.     Write  $\|\cdot \|_s$  for the norm of the weighted Bergman space  
$L^2_a({\bold B},(1-|z|^2)^sdv(z))$.   From the above identity we obtain
$$
\|R^m(f^th)\|_s  \leq   \sum _{k=0}^m|\rho _k|\|f^{t-k}\|_\infty \|R^m(f^kh)\|_s  \leq   \sum _{k=0}^m|\rho _k|\|f^{t-k}\|_\infty \|f^kh\|_{m,s}.
$$
\noindent
Combining the above with the obvious inequality       
$$
\|f^th\|_{m,s} \leq  |f^t(0)h(0)| + \|R^m(f^th)\|_s  \leq   |f(0)|^t\|h\|_{m,s} + \|R^m(f^th)\|_s,
$$
\noindent
we conclude that  $f^t$  is a multiplier for  ${\Cal H}_{m,s}$.     This completes the proof.   $\square $
\medskip
\noindent
{\it Proof of Theorem} 1.4.     Let  $f$  be a multiplier for  ${\Cal H}_{m,s}$, and suppose that  $f$  
does not vanish on  ${\bold B}$.    Thus  $\log f$  exists as an analytic function on  ${\bold B}$.
\medskip
Suppose that   $\|\log f\|_\infty < \infty $.  This  implies that there is a  $c > 0$  such that  $|f(z)| \geq c$  
for every  $z \in {\bold B}$.   Applying Proposition 1.6   to  $f$  and to each  $h \in {\Cal H}_{m,s}$,  we obtain  
$$
R^m((\log f)h)  =  (\log f)R^mh + \sum _{k=0}^ma_kf^{-k}R^m(f^kh),
$$
\noindent
where   $a_0 = -a_1-\cdots -a_m$.  Again, write  $\|\cdot \|_s$  for the norm of the weighted Bergman space  
$L^2_a({\bold B},(1-|z|^2)^sdv(z))$.   From the above identity we obtain
$$
\align
\|R^m((\log f)h)\|_s  &\leq  \|\log f\|_\infty \|R^mh\|_s + \sum _{k=0}^m|a_k|c^{-k}\|R^m(f^kh)\|_s   \\
&\leq   \|\log f\|_\infty \|h\|_{m,s} + \sum _{k=0}^m|a_k|c^{-k}\|f^kh\|_{m,s}.
\endalign
$$
\noindent
Combining the above with the obvious inequality       
$$
\|(\log f)h\|_{m,s} \leq  |(\log f(0))h(0)| + \|R^m((\log f)h)\|_s  \leq   |\log f(0)|\|h\|_{m,s} + \|R^m((\log f)h)\|_s,
$$
\noindent
we conclude that  $\log f$  is a multiplier for  ${\Cal H}_{m,s}$.     
\medskip
On the other hand,  if  $\log f$  is a multiplier for   ${\Cal H}_{m,s}$,  then by Lemma 2.1 we have  
$\|\log f\|_\infty $  $<$  $\infty $.   This completes the proof.   $\square $

\bigskip
\centerline{\bf  3.  Differentiation formulas}    
\medskip
We now turn to the proofs of  Propositions 1.5 and 1.6.
Our approach to the proofs of  Propositions 1.5 and 1.6   is quite different from the approach to the proof of 
(1.2) in [2].                                                                                                                                              
\medskip
Our main idea for the proofs of Propositions 1.5 and 1.6 is to take advantage of the famous 
Faa di Bruno formula.  But instead of the actual formula commonly attributed to Faa di Bruno (see, e.g., [8]),  
where all the coefficients are explicit, all we need is a very crude version of it.  In fact, this is a point 
that we want to emphasize: the proofs of Propositions 1.5 and 1.6 do not need the delicate combinatorics that makes 
the coefficients in the Faa di Bruno formula explicit.  We will see that the proofs of Propositions 1.5 and 1.6  
involve nothing but linear algebra!
\medskip
Given any  $\nu \in {\bold N}$,  we define   $A_\nu $  to be the collection of   $\alpha = (\alpha _1,\dots ,\alpha _\nu ) \in {\bold Z}_+^\nu $  
satisfying the condition
$$
1\cdot \alpha _1 + 2\cdot \alpha _2 + \cdots + \nu \cdot \alpha _\nu   =  \nu .
$$
\noindent  
We follow the usual multi-index notation (see, e.g., [11, page 3]): 
for  $\alpha = (\alpha _1,\dots ,\alpha _\nu )$  $\in $  ${\bold Z}_+^\nu $,  
we write  $|\alpha | = \alpha_1 + \cdots +\alpha _\nu $  and  $(z_1,\dots ,z_\nu )^\alpha = z_1^{\alpha _1}\cdots z_\nu ^{\alpha _\nu }$. 
\medskip
For each  $\nu \in {\bold N}$,  there exists a set of coefficients   $\{b_\alpha : \alpha \in A_\nu \}$  such that the identity 
$$
{d^\nu \over dx^\nu }F(G(x))  =   \sum _{\alpha \in A_\nu }b_\alpha F^{(|\alpha |)}(G(x))(G^{(1)}(x),\dots ,G^{(\nu )}(x))^\alpha 
\tag 3.1
$$
\noindent
holds for all functions  $F$, $G$  that are sufficiently smooth.
\medskip
In contrast to the coefficient-explicit version of the Faa di Bruno formula, 
(3.1) follows easily from an induction on  $\nu $.   For the proofs of 
of Propositions 1.5 and 1.6, the actual values of the  coefficients   $\{b_\alpha : \alpha \in A_\nu \}$  are not relevant; 
the only thing that is relevant is the fact that such coefficients do exist.
\medskip
To prove Propositions 1.5 and 1.6,  for every pair of   $1 \leq i,k \leq m$,  we define
$$
\beta _{i,k}  =  
\left\{
\matrix
k!/(k-i)!  &\text{if}   &i \leq k  \\
\ \ \\
0  &\text{if}   &i > k 
\endmatrix
\right. .
\tag 3.2
$$
\noindent
Then the  $m\times m$  upper-triangular matrix  $[\beta _{i,k}]_{i,k=1}^m$  is invertible.  Therefore there is an  
$m\times m$  matrix   $[c_{k,r}]_{k,r=1}^m$  of rational entries  such that for all   $1 \leq i, r \leq m$,
$$
\sum _{k=1}^m\beta _{i,k}c_{k,r}  =  
\left\{
\matrix
1  &\text{if}   &i = r  \\ 
\ \ \\
0  &\text{if}   &i \neq r
\endmatrix
\right.  . 
\tag 3.3
$$
\noindent
Define the set  $A_{\nu ,i} = \{\alpha \in A_\nu : |\alpha | = i\}$  
for each pair of  $i \leq \nu $   in  $\{1,\dots ,m\}$. 
We write  $C^m_\nu $  for the binomial coefficient  ${m!\over \nu !(m-\nu )!}$.
\medskip
\noindent
{\bf Lemma 3.1.}      {\it Let}  $f$, $h$  {\it be analytic on}  ${\bold B}$.  
{\it If}  $f$  {\it does not vanish on}  ${\bold B}$,  {\it we define}
$$
X_{f,i}h  =  f^{-i}\sum _{\nu = i}^mC_\nu ^m
\sum _{\alpha \in A_{\nu ,i}}b_\alpha (Rf,\dots ,R^\nu f)^\alpha R^{m-\nu }h
$$
\noindent
{\it for}   $1 \leq i \leq m$.   {\it Then  for every}  $1 \leq r \leq m$,  {\it we have}
$$
\sum _{k=1}^m{c_{k,r}\over f^{k}}\sum _{\nu = 1}^mC_\nu ^m(R^\nu f^k)(R^{m-\nu }h)    
= X_{f,r}h.
\tag 3.4
$$
\medskip
\noindent
{\it Proof}.   For each  $1 \leq \nu \leq m$,  if we apply (3.1) to the function  $F(x) = x^k$,  we obtain
$$
\align
R^\nu f^k  &=  \sum _{i=1}^{\min \{k,\nu \}}\sum _{\alpha \in A_{\nu ,i}}b_\alpha {k!\over (k-|\alpha |)!}f^{k-|\alpha |}(Rf,\dots ,R^\nu f)^\alpha  \\
&=  \sum _{i=1}^\nu \beta _{i,k}f^{k-i}\sum _{\alpha \in A_{\nu ,i}}b_\alpha (Rf,\dots ,R^\nu f)^\alpha .
\endalign
$$
\noindent
Thus
$$
\align
\sum _{\nu = 1}^mC_\nu ^m(R^\nu f^k)&(R^{m-\nu }h)    
=  \sum _{\nu = 1}^mC_\nu ^m \sum _{i=1}^\nu \beta _{i,k}f^{k-i}\sum _{\alpha \in A_{\nu ,i}}b_\alpha (Rf,\dots ,R^\nu f)^\alpha R^{m-\nu }h  \\
&=  f^k\sum _{i=1}^m\beta _{i,k}f^{-i}\sum _{\nu = i}^mC_\nu ^m
\sum _{\alpha \in A_{\nu ,i}}b_\alpha (Rf,\dots ,R^\nu f)^\alpha R^{m-\nu }h  \\
&=    f^k\sum _{i=1}^m\beta _{i,k}X_{f,i}h.
\endalign
$$
\noindent
Hence  for each  $1 \leq r \leq m$,
$$
\sum _{k=1}^m{c_{k,r}\over f^k}\sum _{\nu = 1}^mC_\nu ^m(R^\nu f^k)(R^{m-\nu }h)  
=  \sum _{i=1}^m\sum _{k=1}^mc_{k,r}\beta _{i,k}X_{f,i}h = X_{f,r}h,  
$$
\noindent
where the second  $=$  is obtained from  (3.3).    $\square $

\medskip
\noindent
{\it Proof of Proposition} 1.5. 
For each natural number  $1 \leq r \leq m$,  we define
$$
u(r) = \prod _{j=0}^{r-1}(t-j) = (-1)^r\prod _{j=0}^{r-1}(j-t).
\tag 3.5
$$
\noindent
Thus   $(x^t)^{(r)} = u(r)x^{t-r}$,  $1 \leq r \leq m$.
We then define  $\rho _k = \sum _{r=1}^mu(r)c_{k,r}$   for each  $1 \leq k \leq m$.
Let us verify that (1.3) holds for the   $\rho _1, \dots ,\rho _m$  so defined.
\medskip
Let   $f, h$ be analytic functions on  ${\bold B}$,  and suppose that  $f$  does not vanish on  ${\bold B}$.  
By the Leibniz rule for differentiation,  we have
$$
\align
R^m(f^th) - f^tR^mh    &=    \sum _{\nu = 1}^mC_\nu ^m(R^\nu f^t)(R^{m-\nu }h)   \quad  \text{and}  \\
R^m(f^kh) - f^kR^mh  &=    \sum _{\nu = 1}^mC_\nu ^m(R^\nu f^k)(R^{m-\nu }h),  \quad  1 \leq k \leq m.
\endalign
$$
Thus (1.3) will hold  for the  $\rho _1, \dots ,\rho _m$  defined above if we can show that
$$
\sum _{\nu = 1}^mC_\nu ^m(R^\nu f^t)(R^{m-\nu }h) 
=  f^t\sum _{r=1}^mu(r)\sum _{k=1}^m{c_{k,r}\over f^k}\sum _{\nu = 1}^mC_\nu ^m(R^\nu f^{k})(R^{m-\nu }h).
\tag 3.6 
$$   
\medskip
To prove (3.6),  we apply (3.1)  to the function  $F(x) = x^t$.  For each  
$1 \leq \nu \leq m$,  (3.1)  gives us  
$$
R^\nu f^t  =  \sum _{\alpha \in A_\nu }b_\alpha u(|\alpha |)f^{t-|\alpha |}(Rf,\dots ,R^\nu f)^\alpha 
$$
\noindent
(see (3.5)).  Therefore  
$$
\align
\sum _{\nu = 1}^mC_\nu ^m(R^\nu f^t)&(R^{m-\nu }h)    
=  \sum _{\nu = 1}^mC_\nu ^m \sum _{\alpha \in A_\nu }b_\alpha u(|\alpha |)f^{t-|\alpha |}
(Rf,\dots ,R^\nu f)^\alpha R^{m-\nu }h  \\
&=  \sum _{\nu = 1}^mC_\nu ^m \sum _{r=1}^\nu \sum _{\alpha \in A_{\nu ,r}}b_\alpha 
u(|\alpha |)f^{t-|\alpha |}(Rf,\dots ,R^\nu f)^\alpha R^{m-\nu }h  \\
&=  \sum _{r=1}^mu(r)f^{t-r}\sum _{\nu = r}^mC_\nu ^m
\sum _{\alpha \in A_{\nu ,r}}b_\alpha (Rf,\dots ,R^\nu f)^\alpha R^{m-\nu }h  \\
&=    f^t\sum _{r=1}^mu(r)X_{f,r}h. 
\tag 3.7
\endalign
$$
\noindent
Applying Lemma 3.1 in  (3.7),  we obtain (3.6).  This completes the proof.   $\square $
\medskip
\noindent
{\it Proof of Proposition} 1.6. 
For each natural number  $1 \leq r \leq m$,  we define
$$
v(r) = (r-1)!(-1)^{r-1}.
\tag 3.8
$$
\noindent
Thus   $(\log x)^{(r)} = v(r)x^{-r}$,  $1 \leq r \leq m$.
We then define  $a_k = \sum _{r=1}^mv(r)c_{k,r}$   for each  $1 \leq k \leq m$.
Let us verify that (1.4) holds for the   $a_1, \dots ,a_m$  so defined.
\medskip
Let   $f, h$ be analytic functions on  ${\bold B}$,  and suppose that  $f$  does not vanish on  ${\bold B}$.  
By the Leibniz rule for differentiation,  we have
$$
\align
R^m((\log f)h) - (\log f)R^mh    &=    \sum _{\nu = 1}^mC_\nu ^m(R^\nu (\log f))(R^{m-\nu }h)   \quad  \text{and}  \\
R^m(f^kh) - f^kR^mh  &=    \sum _{\nu = 1}^mC_\nu ^m(R^\nu f^k)(R^{m-\nu }h),  \quad  1 \leq k \leq m.
\endalign
$$
\noindent
Thus (1.4) will hold  for the  $a_1, \dots ,a_m$  defined above if we can show that
$$
\sum _{\nu = 1}^mC_\nu ^m(R^\nu (\log f))(R^{m-\nu }h) 
=  \sum _{r=1}^mv(r)\sum _{k=1}^m{c_{k,r}\over f^k}\sum _{\nu = 1}^mC_\nu ^m(R^\nu f^{k})(R^{m-\nu }h).
\tag 3.9 
$$   
\medskip
To prove (3.9),  we apply (3.1)  to the function  $F(x) = \log x$.  For each  
$1 \leq \nu \leq m$,  (3.1)  gives us  
$$
R^\nu (\log f)  =  \sum _{\alpha \in A_\nu }b_\alpha v(|\alpha |)f^{-|\alpha |}(Rf,\dots ,R^\nu f)^\alpha 
$$
\noindent
(see (3.8)).  Therefore  
$$
\align
\sum _{\nu = 1}^mC_\nu ^m(R^\nu (\log f))&(R^{m-\nu }h)    
=  \sum _{\nu = 1}^mC_\nu ^m \sum _{\alpha \in A_\nu }b_\alpha v(|\alpha |)f^{-|\alpha |}
(Rf,\dots ,R^\nu f)^\alpha R^{m-\nu }h  \\
&=  \sum _{\nu = 1}^mC_\nu ^m \sum _{r=1}^\nu \sum _{\alpha \in A_{\nu ,r}}b_\alpha 
v(|\alpha |)f^{-|\alpha |}(Rf,\dots ,R^\nu f)^\alpha R^{m-\nu }h  \\
&=  \sum _{r=1}^mv(r)f^{-r}\sum _{\nu = r}^mC_\nu ^m
\sum _{\alpha \in A_{\nu ,r}}b_\alpha (Rf,\dots ,R^\nu f)^\alpha R^{m-\nu }h  \\
&=   \sum _{r=1}^mv(r)X_{f,r}h. 
\tag 3.10
\endalign
$$
\noindent
Applying Lemma 3.1 in  (3.10),  we obtain (3.9).  This completes the proof.   $\square $

\bigskip
\centerline{\bf  4.  More explicit formulas for the coefficients}
\medskip
Note that Propositions 1.5 and 1.6 were proved without any kind of computation.  We will now add up the sums    
$\rho _k = \sum _{r=1}^mu(r)c_{k,r}$  and  $a_k = \sum _{r=1}^mv(r)c_{k,r}$  to make these coefficients more 
explicit, which involves some elementary calculation.
\medskip
First, note that   $[\beta _{i,k}] = e^NF$,  where  $F$  is the diagonal matrix whose diagonal entries are  
$1!, \dots ,m!$  and   $N = [s_{i,k}]$,  where  $s_{i,i+1} = 1$  for   $1 \leq i \leq m-1$  and  $s_{i,k} = 0$  if   $k \neq i+1$.   
Consequently,  $[c_{k,r}] = F^{-1}e^{-N}$.  That is, for all  $k, r \in \{1,\dots ,m\}$  we have
$$
c_{k,r}  =  
\left\{
\matrix
{(-1)^{r-k}\over k!(r-k)!}    &\text{if}     &r \geq k  \\
\  \  \\
0   &\text{if}     &r < k 
\endmatrix
\right. .
\tag 4.1
$$
\noindent
Next, recall that for all  non-negative integers $p \geq 0$  and  $q \geq 0$,  there is the identity
$$
\sum _{j=0}^q{(j+p)!\over j!}  =  {(q+p+1)!\over q!(p+1)},
\tag 4.2
$$
\noindent
which is proved by an easy induction on  $q$.
\medskip
By (3.8), (4.1) and (4.2),  for each  $1 \leq k \leq m$   we have
$$
\align
a_k &= \sum _{r=k}^{m}(-1)^{r-1+r-k}{(r-1)!\over k!(r-k)!}  =  {(-1)^{k+1}\over k!}\sum _{j=0}^{m-k}{(j+k-1)!\over j!}   \\
&=  {(-1)^{k+1}\over k!}\cdot {(m-k+k-1+1)!\over (m-k)!k}   
=  {(-1)^{k+1}\over k}\cdot {m!\over k!(m-k)!}.
\endalign
$$
\noindent
Similarly,  by (4.1) and (3.5),  for each  $1 \leq k \leq m$,
$$
\rho _k = \sum _{r=k}^m(-1)^{r+r-k}{\prod _{i=0}^{r-1}(i-t)\over k!(r-k)!}  
=  {(-1)^k\over k!}\sum _{j=0}^{m-k}{\prod _{i=0}^{j+k-1}(i-t)\over j!}.
\tag 4.3
$$
\noindent
In the case where $t$  is a negative integer,  by (4.2) we have
$$
\align
\sum _{j=0}^{m-k}{\prod _{i=0}^{j+k-1}(i-t)\over j!}  &=  \sum _{j=0}^{m-k}{(j+k-1-t)!\over (-t-1)!j!}   
=  {(m-k+k-1-t +1)!\over (-t-1)!(m-k)!(k-1-t+1)}  \\
&=  {(m-t)!\over (-t-1)!(m-k)!(k-t)}   =  {1\over (m-k)!}\prod _{i\in \{0,1,\dots,m\}\backslash \{k\}}(i - t).
\tag 4.4
\endalign
$$
\noindent
If  $a$  and  $b$  are polynomials such that  $a(t) = b(t)$  for every 
negative integer  $t$,  then   $a(t) = b(t)$   for every  $t \in {\bold R}$.  Thus it follows from (4.3) and (4.4) that  
for every  $1 \leq k \leq m$,
$$
\rho _k   =  {(-1)^k\over k!(m-k)!}\prod _{i\in \{0,1,\dots,m\}\backslash \{k\}}(i - t).
$$

\bigskip
\centerline{\bf   References}
\medskip
\noindent
1.   W. Arveson,  Subalgebras of $C^\ast $-algebras. III.  Multivariable operator theory,  Acta Math. {\bf 181}  (1998), 159-228.

\noindent
2.  G. Cao, L. He and K. Zhu,  Spectral theory of multiplication operators on Hardy-Sobolev spaces, J. Funct. Anal. 
{\bf 275}  (2018),  1259-1279.

\noindent
3. S.  Costea,  E. Sawyer and B. Wick,  The corona theorem for the Drury-Arveson Hardy space and other holomorphic 
Besov-Sobolev spaces on the unit ball in ${\bold C}^n$,   Anal. PDE {\bf 4} (2011),   499-550.

\noindent
4.  S. Drury, A generalization of von Neumann's inequality to the complex ball, 
Proc. Amer. Math. Soc. {\bf 68} (1978), 300-304.

\noindent
5.  Q. Fang and J. Xia,  Corrigendum to ``Multipliers and essential norm on the Drury-Arveson space", 
Proc. Amer. Math. Soc. {\bf 141} (2013), 363-368. 

\noindent
6.  Q. Fang and J. Xia,  A hierarchy of von Neumann inequalities?, J. Operator Theory {\bf 72} (2014), 219-239.

\noindent
7.  Q. Fang and J. Xia,  Analytical aspects of the Drury-Arveson space, Handbook of analytic operator theory, 203-221, 
CRC Press/Chapman Hall Handb. Math. Ser., CRC Press, Boca Raton, FL, 2019.

\noindent
8.  W. Johnson, The curious history of Fa\`a di Bruno's formula, Amer. Math. Monthly {\bf 109} (2002),  217-234.

\noindent
9.  A.  Lubin,  Weighted shifts and products of subnormal operators,  
Indiana Univ. Math. J. {\bf 26} (1977),   839-845.

\noindent
10.  S.  Richter  and  J.  Sunkes,   Hankel operators, invariant subspaces, and cyclic vectors in the Drury-Arveson space,  
Proc. Amer. Math. Soc. {\bf 144} (2016),  2575-2586. 

\noindent
11.  W.  Rudin, Function theory in the unit ball of ${\text{\bf C}}^n$,  Springer-Verlag, New York, 1980.

\noindent
12.  K. Zhu,  Spaces of holomorphic functions in the unit ball, Graduate Texts in Mathematics {\bf 226},  Springer-Verlag, New York, 2005.

\bigskip
\bigskip
\bigskip
\bigskip
\noindent
Jingbo Xia

\noindent
College of Data Science, Jiaxing University,  Jiaxing 314001,  China

\noindent
and

\noindent
Department of Mathematics, State University of New York at Buffalo, Buffalo, NY 14260,  USA  

\noindent
E-mail: \url{jxia@acsu.buffalo.edu}

\bigskip
\noindent
Congquan Yan

\noindent
College of Data Science, Jiaxing University,  Jiaxing 314001,  China

\noindent
E-mail: \url{yancongquan@zjxu.edu.cn}

\bigskip
\noindent
Danjun Zhao

\noindent
College of Data Science, Jiaxing University,  Jiaxing 314001,  China

\noindent
E-mail: \url{danjunzhao@zjxu.edu.cn}

\bigskip
\noindent
Jingming Zhu

\noindent
College of Data Science, Jiaxing University,  Jiaxing 314001,  China

\noindent
E-mail: \url{jingmingzhu@zjxu.edu.cn}

\end